\documentclass{article}
\usepackage{amsmath, amscd, amsfonts, amssymb, amsthm}

\newtheorem{theorem}{Theorem}[section]
\newtheorem{corollary}[theorem]{Corollary}
\newtheorem{lemma}[theorem]{Lemma}
\begin{document}

\title{The Critical Value of the Contact Process with Added and Removed Edges}
\author{Paul Jung}

\date{\textit{Department of Mathematics, Cornell University}} \maketitle

\begin{abstract}
We show that the critical value for the contact process on a
vertex-transitive graph $\mathcal{G}$ with finitely many edges
added and/or removed is the same as the critical value for the
contact process on $\mathcal{G}$.  This gives a partial answer to
a conjecture of Pemantle and Stacey.
\end{abstract}


\section{Introduction}
The contact process with infection rate $\lambda$ is an
interacting particle system whose state space is
$X=\{0,1\}^{\mathcal{G}}$ for some connected graph $\mathcal{G}$
with countably many vertices (we use $\mathcal{G}$ to describe
both the vertex set and the graph, depending on the context). This
process was first introduced by Harris (1974) as a stochastic
model for the spread of an infection (1's mark infected sites and
0's mark healthy sites). A simple description of the process is as
follows. Each healthy site $x\in\mathcal{G}$ becomes infected at
an exponential rate which is proportional to the number of
infected neighbors of $x$ and each infected site becomes healthy
at exponential rate $1$. The contact process has been extensively
studied since the 1970's; for a detailed account of what is known
about this process we refer the reader to Liggett(1999).

To define the process more precisely, let $\eta\in X$. We denote
\begin{equation*}
\eta_{x}(u)=\left\{
\begin{array}{ll}
1-\eta(u)&\text{if }u=x\\
\eta(u)&\text{if }u\neq x\\
\end{array}
\right.
\end{equation*}
and
\begin{equation*}
c_\mathcal{G}(x,\eta)=\left\{
\begin{array}{ll}
1&\text{if }\eta(x)=1\\
\lambda\sum_{|y-x|=1}{\eta(y)}&\text{if }\eta(x)=0
\end{array}
\right..
\end{equation*}
The generator of the contact process $\eta_t$ is now given by the
closure of the operator $L_{\mathcal{G}}$ on $\mathcal{D}(X)$, the
set of all functions on $X$ depending on finitely many
coordinates:
\begin{equation*}\label{generator}
L_{\mathcal{G}}
f(\eta)=\sum_{x}{c_\mathcal{G}(x,\eta)[f(\eta_{x})-f(\eta)]},\,
f\in\mathcal{D}(X).
\end{equation*}
The subscript $\mathcal{G}$ reminds us that the generator and the
transition rates depend on the graph $\mathcal{G}$.  We write
$S_{\mathcal{G}}(t)$ for the semigroup of this process.

If the initial state of the process is such that $\eta_0(x)=1$ for
only finitely many $x\in\mathcal{G}$, then the process is just a
Markov chain, and its state at time $t\ge 0$ is given by $A_t$
where $A$ is the finite subset of $\mathcal{G}$ which is exactly
the set of all infected sites at time $t$. In the sequel, $A$ will
always be a finite subset of $\mathcal{G}$. One of the basic tools
used in the study of the contact process is its self-duality.  In
other words:
\begin{equation}\label{dual}
P^{\eta}[\eta_t(x)=0\text{ for all }x\in A]=P^A[\eta(x)=0\text{
for all }x\in A_t].
\end{equation}
The above result is based on a graphical construction of the
contact process, and it can be found in any reference on the
contact process (for example Liggett(1999)).

One of the goals in the treatment of any interacting particle
system is to study the invariant measures.  For the contact
process, it is easily seen that $\delta_0$, the point mass on all
$0$'s, is an invariant measure. Also, if $\delta_1$ is the point
mass on all $1$'s then simple coupling and monotonicity arguments
show that the \textit{upper invariant measure}
$$\lim_{t\rightarrow\infty}\delta_1
S_{\mathcal{G}}(t)={\mu}_{\mathcal{G}}$$ exists and
${\mu}_{\mathcal{G}}$ stochastically dominates every other
invariant measure. If ${\mu}_{\mathcal{G}}=\delta_0$ then
$\delta_0$ is the only invariant measure and the process is
ergodic. We have the following definition for $\lambda_c$, the
critical value of the contact process:
$$\lambda_c=\sup\{\lambda\ge 0: \mu_\mathcal{G}=\delta_0\}.$$
We note here that $\lambda_c$ is often called the \textit{global
survival} critical value or the \textit{lower} critical value. We
mention this because there is another natural critical value for
the contact process which is known as the \textit{local survival}
critical value or the \textit{upper} critical value, however, in
this paper we have no need to define this other critical value.

Suppose $\mathcal{G}$ is a connected graph with countably many
vertices.  Let $\mathcal{G}'$ be a graph formed by adding $n$
edges to the graph $\mathcal{G}$.  We will say that the
$i^\text{th}$ new edge is placed between the vertices $u_i$ and
$v_i$ where $u_i\neq v_i$ since loops are meaningless in the
contact process. The new edges can be placed between two vertices
that already have an edge in $\mathcal{G}$, however, to simplify
things we assume that each element in the set $\{u_1, v_1,\ldots,
u_n, v_n\}$ is distinct. We use the phrase ``to simplify things"
here because this requirement is not necessary, but it makes the
proofs easier to follow.

Let $\lambda_c'$ be the critical value of the contact process on
$\mathcal{G}'$. If $\lambda_c$ is the critical value for the
contact process on $\mathcal{G}$, then Pemantle and Stacey(2000)
have conjectured that $\lambda_c'=\lambda_c$. When
$\mathcal{G}=\mathbb{Z}^d$, the following argument which uses
duality together with a result of Bezuidenhout and Grimmett(1991)
shows that this is true.  For simplicity we consider the case
where $\mathcal{G}'$ differs from $\mathcal{G}$ only by the
addition of one edge between the vertices $u$ and $v$.

Let $A_t$ be the process on $\mathcal{G}$ and let $A_t'$ be the
process on $\mathcal{G}'$. Now suppose there exists a nontrivial
upper invariant measure $\mu_{\mathcal{G}'}$ for some
$\lambda<\lambda_c$. If $||\cdot ||$ denotes graph distance from
some distinguished vertex labelled the origin, then
\begin{eqnarray*}
\mu_{\mathcal{G}'}\{\eta:\eta(x)=1\}&=&{P}^{\{x\}}_\lambda(\{A_t'\neq
\emptyset\text{ for all
}t\})\\
&=&{P}^{\{x\}}_\lambda(\{A_t'\neq\emptyset\text{ for all }t\}\cap
\{u\in A_t' \text{ or }v\in A_t'\text{ for some }t\})\\
&\le&e^{-c||x||}\text{ for some }c>0
\end{eqnarray*}
where the inequality comes from Theorem 1.7 of Bezuidenhout and
Grimmett(1991).  But now the above inequality implies that
$\sum_x\mu_{\mathcal{G}'}\{\eta:\eta(x)=1\}<\infty$ which means
that $\mu_{\mathcal{G}'}$ concentrates on configurations with
finitely many ones contradicting its invariance.  Therefore
$\lambda_c'=\lambda_c$.

The above argument can easily be extended to adding any finite
number of edges to $\mathcal{G}=\mathbb{Z}^d$.  The issue of
removing edges is a bit trickier; it requires us to start from a
graph $\mathcal{G}$ which is exactly $\mathbb{Z}^d$ with finitely
many edges removed. This is not much of a problem since it seems
that the argument used to show exponential decay in Theorem 1.7 in
Bezuidenhout and Grimmett(1991) can be extended to such graphs.
However, we would also like to know that the conjecture of
Pemantle and Stacey holds for graphs such as $\mathbb{T}^d$, and
here we run into a problem since the above argument depends on the
amenability of $\mathbb{Z}^d$. The goal of this paper is to
introduce an alternate argument which shows that
$\lambda_c'=\lambda_c$ whenever we can satisfy a certain
integrability condition known to hold for the subcritical contact
process even on nonamenable graphs.

Since $\mathcal{G}$ and $\mathcal{G}'$ have the same vertex sets,
we will always use the notation $\mathcal{G}$ when referring to
the vertex set of either graph. As noted above, it should be clear
from the context when $\mathcal{G}$ refers to the vertex set
rather than the graph.

\begin{theorem}\label{theorem}
If $A_t^{\{o\}}$ is the contact process on $\mathcal{G}$ starting
from one infection at the origin, $o$, and for $\lambda<\lambda_c$
we have
\begin{equation}\label{contact3}
\int_0^\infty E|A_t^{\{o\}}| dt<\infty,
\end{equation}
then
$\lambda_c'=\lambda_c$.
\end{theorem}

Note that $\int_0^\infty E|A_t^{\{o\}}| dt<\infty$ if and only if
$\int_0^\infty E|A_t^{A}| dt<\infty$ for all finite
$A\subset\mathcal{G}$.

\begin{corollary}\label{cor}
If $\mathcal{G}$ is a vertex-transitive graph and $\mathcal{G}'$
is formed by adding and/or removing finitely many edges from
$\mathcal{G}$, then the critical values for the contact process on
$\mathcal{G}$ and $\mathcal{G}'$ are the same.
\end{corollary}

The corollary follows from the arguments of Aizenman and
Barsky(1987) where it is shown that (\ref{contact3}) holds for all
$\lambda<\lambda_c$ whenever $\mathcal{G}$ is transitive. Aizenman
and Barsky(1987) actually concentrate on a discrete-time
percolation model, but we have learned through personal
communication that Aizenman has an unpublished extension to the
contact process.  The current author also has proved such an
extension in a forthcoming paper. It can be seen in the arguments
of Aizenman and Barsky(1987), that the proof is also valid for
transitive graphs with finitely many edges removed.

Some comments are in order concerning the techniques used in the
proofs below. Theorem \ref{theorem} is a result that gives
information about the contact process when the infection rates are
perturbed at a finite number of sites. The techniques used in the
proof work equally well when the healing rates are perturbed
(these processes are known as inhomogeneous contact processes). In
particular, it can be seen from the proofs below that if the
healing rates are lowered at a finite number of sites for the
contact process on a vertex-transitive graph, then the critical
value is left unchanged.  This is a special case of a result
proved by Madras, Schinazi, and Schonmann(1994) for inhomogeneous
contact processes on $\mathbb{Z}^d$.  We also note that the
techniques in the proofs can also be applied to other spin
systems.  For example, one can extend to all dimensions, Theorem 2
of Handjani(1999) which concerns a perturbed biased-voter model on
$\mathbb{Z}$.

\section{A generator computation}\label{im}
In this section we prove a crucial lemma which uses certain
coupled processes to gain information about the evolution of the
upper invariant measure for the contact process on $\mathcal{G}'$
(which we denote as ${\mu}_{\mathcal{G}'}$) under the semigroup
$S_{\mathcal{G}}$. As we will see, the lemma basically boils down
to a generator computation, giving us the title of this section.
The methods of this section are motivated by the infinitesimal
coupling of the exclusion process used in Andjel, Bramson, and
Liggett(1988) and again in Jung(2004).

Before stating the lemma, we describe the couplings to be used.
All couplings to be used will follow the motion of the basic
coupling for two processes $\eta_t$ and $\xi_t$. The \textit{basic
coupling} is the coupling of $\eta_t$ and $\xi_t$ which allows the
two processes to move together as much as possible (see
Liggett(1985) Chapter III for more details). One of the most
useful properties of the basic coupling is the fact that it
preserves stochastic domination in time. In particular, if
$\eta_t$ and $\xi_t$ are coupled using a basic coupling then
$\eta_0(x)\le\xi_0(x)$ for all $x\in\mathcal{G}$ implies that
$\eta_t(x)\le\xi_t(x)$ for all $x\in\mathcal{G}$.

We now describe the various initial measures we use for the
couplings below. Let $\{u_i\}$ and $\{v_i\}$ be as before.
If $\eta$ is given the measure ${\mu}_{\mathcal{G}'}$, define
\begin{eqnarray*}
D_{u}^i=\{\eta: \eta(u_i)=0,\eta(v_i)=1\}\,\,\text{ and
}\,\,D_v^i=\{\eta: \eta(u_i)=1,\eta(v_i)=0\}.
\end{eqnarray*}
The measures $\mu_{u}^i,$ and $\mu_{v}^i$ are defined by
conditioning $\mu_{\mathcal{G}'}$ on the events $D_u^i$ and
$D_v^i$ respectively. Also, define the measures $\hat{\mu}_{u}^i$
and $\hat{\mu}_{v}^i$ to be exactly equal to ${\mu}_{u}^i$ and
${\mu}_{v}^i$ except that we change the values at $u_i$ and $v_i$
so that $\eta(u_i)=\eta(v_i)=1$.

We can now define the initial measures $\nu_u^i$ and $\nu_v^i$ for
the coupled processes $(\eta_t^{u,i},\xi_t^{u,i})$ and
$(\eta_t^{v,i},\xi_t^{v,i})$. The measure $\nu_z^i$ for $z=u,v$
has marginal measures $\mu_z^i$ and $\hat{\mu}_z^i$ corresponding
to $\eta_0^{z,i}$ and $\xi_0^{z,i}$ respectively and the marginals
are coupled so that $\eta_0^{z,i}(x)\le\xi_0^{z,i}(x)$ for all
$x$.

\begin{lemma}\label{coupling2}
\begin{equation*}
 \frac{d}{dt}\mu_{\mathcal{G}'}
{S}_\mathcal{G}(t)\{\eta:\eta(x)=1\}= \sum_{i=1}^n\sum_{z=u,v}
{\lambda\mu_{\mathcal{G}'}\{D_z^i\}}E
[{\eta}_t^{z,i}(x)-{\xi}_t^{z,i}(x)]
\end{equation*}
\end{lemma}

\begin{proof}
We write $1_x(\eta)=\eta(x)$ and
$1_x^t(\eta)=S_{\mathcal{G}}(t)1_x(\eta)$. Letting
$\mu_{\mathcal{G}'}^t=\mu_{\mathcal{G}'} {S}_{\mathcal{G}}(t)$ we
have
\begin{eqnarray*}
\frac{d}{dt}\mu_{\mathcal{G}'}
{S}_\mathcal{G}(t)\{\eta:\eta(x)=1\} &=&\lim_{s\rightarrow
0}\frac{1}{s}[\int 1_x\,
d\mu_{\mathcal{G}'}^{t+s}-\int 1_x\,d\mu_{\mathcal{G}'}^t]\\
&=&\lim_{s\rightarrow 0}\frac{1}{s}[\int 1^t_x\,
d\mu_{\mathcal{G}'}^{s}-\int 1^t_x\,d\mu_{\mathcal{G}'}].
\end{eqnarray*}
By the definition of the generator, the above is equal to
\begin{eqnarray*}
&=&\int {L_{\mathcal{G}}} 1^t_x \,d\mu_{\mathcal{G}'}\\
&=& \int L_{\mathcal{G}'} 1^t_x \, d\mu_{\mathcal{G}'}
+\sum_{i=1}^n\sum_{z=u,v}\int
(c_{\mathcal{G}}(z_i,\eta)-c_{\mathcal{G}'}(z_i,\eta))[1^t_x(\eta_{z_i})-1^t_x(\eta)]\,
d\mu_{\mathcal{G}'}\\
&=&\sum_{i=1}^n\left[\int
\eta(v_i)(1-\eta(u_i))(-\lambda)[1^t_x(\eta_{u_i})-1^t_x(\eta)]\,
d\mu_{\mathcal{G}'} +\int
\eta(u_i)(1-\eta(v_i))(-\lambda)[1^t_x(\eta_{v_i})-1^t_x(\eta)]\, d\mu_{\mathcal{G}'}\right]\\
&=&\sum_{i=1}^n\lambda\left[\int
\eta(v_i)(1-\eta(u_i))[1^t_x(\eta)-1^t_x(\eta_{u_i})]\,
d\mu_{\mathcal{G}'}
+\int \eta(u_i)(1-\eta(v_i))[1^t_x(\eta)-1^t_x(\eta_{v_i})]\, d\mu_{\mathcal{G}'}\right]\\
&=&\sum_{i=1}^n\sum_{z=u,v}\lambda\mu_{\mathcal{G}'}\{D_z^i\}\int
[1^t_x(\eta)-1^t_x(\eta_{z_i})] \, d\mu_z^i\\
&=&\sum_{i=i}^n\sum_{z=u,v} {\lambda\mu_{\mathcal{G}'}\{D_z^i\}}E
[{\eta}_t^{z,i}(x)-{\xi}_t^{z,i}(x)].
\end{eqnarray*}
The third equality above follows since $\int L_{\mathcal{G}'}
1^t_x \, d\mu_{\mathcal{G}'}=0$ by the invariance of
$\mu_{\mathcal{G}'}$ under $L_{\mathcal{G}'}$.
\end{proof}


\section{Proof of Theorem \ref{theorem}}

\begin{proof}[Proof of Theorem \ref{theorem}]
Let $\zeta_t^{z,i}=\xi_t^{z,i}-\eta_t^{z,i}$. It is clear from the
way that $\eta_t^{z,i}$ and $\xi^{z,i}_t$ are coupled that
$\zeta_t^{z,i}(x)\le A_t^{\{z_i\}}(x)$ for all $x\in\mathcal{G}$.
Therefore
\begin{eqnarray*}
\sum_{x\in\mathcal{G}}E[\xi_t^{z,i}(x)-\eta_t^{z,i}(x)]\le
E|A_t^{\{z_i\}}|.
\end{eqnarray*}
By Lemma \ref{coupling2} we get that
\begin{eqnarray}\label{c3}
(-1)\sum_{x\in\mathcal{G}}\frac{d}{dt}\mu_{\mathcal{G}'}
{S}_\mathcal{G}(t)\{\eta(x)=1\}&\le&\sum_{i=1}^n\sum_{z=u,v}
{\lambda\mu_{\mathcal{G}'}\{D_z^i\}}E|A_t^{\{z_i\}}|.
\end{eqnarray}

Using monotonicity arguments it is easy to show that
$\lim_{t\rightarrow\infty}\mu_{\mathcal{G}'}S_{\mathcal{G}}=\mu_\mathcal{G}$
so integrating both sides of (\ref{c3}) with respect to $t$ from
$0$ to $\infty$ gives
\begin{eqnarray*}
\sum_{x\in\mathcal{G}}[\mu_{\mathcal{G}'}\{\eta(x)=1\}-\mu_{\mathcal{G}}\{\eta(x)=1\}]&=&
(-1)\int_0^\infty\sum_{x\in\mathcal{G}}\frac{d}{dt}\mu_{\mathcal{G}'}
{S}_\mathcal{G}(t)\{\eta(x)=1\}\, dt\\
&\le& \int_0^\infty \sum_{i=1}^n\sum_{z=u,v}
{\lambda\mu_{\mathcal{G}'}\{D_z^i\}}E|A_t^{\{z_i\}}|\,dt
\end{eqnarray*}
When $\lambda<\lambda_c$, the right-hand side is finite and
$\mu_\mathcal{G}=\delta_0$, therefore
$$\sum_{x\in\mathcal{G}}\mu_{\mathcal{G}'}\{\eta(x)=1\}<\infty.$$
But this implies that when $\lambda<\lambda_c$,
$\mu_{\mathcal{G}'}$ concentrates on configurations with finitely
many infected sites. Since $\mu_{\mathcal{G}'}$ is a stationary
distribution for a Markov chain which has $\delta_0$ as its only
absorbing state, it must be that $\mu_{\mathcal{G}'}=\delta_0$
when $\lambda<\lambda_c$ which implies $\lambda_c'=\lambda_c$.
\end{proof}

\textbf{Acknowledgement}. We thank Rick Durrett for many useful
discussions concerning the contact process and for his mentorship
this past year.

\end{document}